\documentclass[11pt]{amsart}
\usepackage{ amsmath, amsthm, amsfonts, hyperref, graphicx, ifpdf}
\usepackage{mathrsfs,epsfig}
\usepackage[dvipsnames,usenames]{color}
\hypersetup{
   unicode=false,          % non-Latin characters in Acrobat bookmarks
   pdftoolbar=true,        % show Acrobat toolbar?
   pdfmenubar=true,        % show Acrobat menu?
   pdffitwindow=false,     % window fit to page when opened
   pdfstartview={},    % fits the width of the page to the window
   pdftitle={},    % title
   pdfauthor={},     % author
   pdfsubject={},   % subject of the document
   pdfcreator={},   % creator of the document
   pdfproducer={}, % producer of the document
   pdfkeywords={}, % list of keywords
   pdfnewwindow=true,      % links in new window
   colorlinks=true,       % false: boxed links; true: colored links
   linkcolor=blue,          % color of internal links
   citecolor=blue,        % color of links to bibliography
   filecolor=blue,      % color of file links
   urlcolor=blue          % color of external links
}

\newtheorem{theorem}{Theorem}[section]
\newtheorem{cor}[theorem]{Corollary}
\newtheorem{lem}[theorem]{Lemma}
\newtheorem{pro}[theorem]{Proposition}
\newtheorem{remark}[theorem]{Remark}
\newtheorem{Def}[theorem]{Definition}
\theoremstyle{definition}

\newcommand{\af}{angle function}
\newcommand{\ap}{avoidance principle}
\newcommand{\bnb}{\overline\nabla}
\newcommand{\cmc}{constant mean curvature}

\newcommand{\ee}{evolution equation}
\newcommand{\ehs}{evolving hypersurface}
\newcommand{\F}{\frac}
\newcommand{\Fm}{Fuchsian manifold}

\newcommand{\htm}{hyperbolic three-manifold}

\newcommand{\ihs}{initial hypersurface}
\newcommand{\ivp}{initial value problem}
\newcommand{\mc}{mean curvature}
\newcommand{\mcf}{mean curvature flow}

\newcommand{\nb}{\nabla}

\newcommand{\pc}{principal curvature}

\newcommand{\ricci}{Ricci curvature}

\newcommand{\sff}{second fundamental form}
\newcommand{\tg}{totally geodesic}

\newcommand{\wf}{warping function}
\newcommand{\wpp}{warped product}
\newcommand{\wppm}{warped product manifold}
\newcommand{\wrt}{with respect to}

\newcommand{\hg}{hyperbolic geometry}

\newcommand{\be}{\begin{equation}}
\newcommand{\ene}{\end{equation}}
\newcommand{\br}{\begin{remark}}
\newcommand{\er}{\end{remark}}
\newcommand{\bl}{\begin{lem}}
\newcommand{\el}{\end{lem}}
\newcommand{\bcor}{\begin{cor}}
\newcommand{\ecor}{\end{cor}}
\newcommand{\bpro}{\begin{pro}}
\newcommand{\epro}{\end{pro}}
\newcommand{\ben}{\begin{enumerate}}
\newcommand{\een}{\end{enumerate}}
\newcommand{\bp}{\begin{proof}}
\newcommand{\ep}{\end{proof}}
\newcommand{\bpo}{\begin{pro}}
\newcommand{\epo}{\end{pro}}
\newcommand{\beq}{\begin{equation*}}
\newcommand{\eeq}{\end{equation*}}
\newcommand{\bear}{\begin{eqnarray}}
\newcommand{\eear}{\end{eqnarray}}
\newcommand{\beqar}{\begin{eqnarray*}}
\newcommand{\eeqar}{\end{eqnarray*}}
\newcommand{\bt}{\begin{theorem}}
\newcommand{\et}{\end{theorem}}

\newcommand{\nablabar}{{\overline{\nabla}}}

\newcommand{\n}{\mathbf{n}}

\newcommand{\ta}{\Theta}
\newcommand{\vv}{\vec{v}}
\newcommand{\vnu}{\boldsymbol{\nu}}

\newcommand{\R}{\mathbb{R}}

\newcommand{\inner}[2]{\langle #1\,,#2\rangle}
\newcommand{\ddl}[2]{\frac{d{#1}}{d{#2}}}

\newcommand{\ppl}[2]{\frac{\partial{#1}}{\partial{#2}}}

\numberwithin{equation}{section}

\allowdisplaybreaks

\def\XXint#1#2#3{{\setbox0=\hbox{$#1{#2#3}{\int}$}
    \vcenter{\hbox{$#2#3$}}\kern-.5\wd0}}

\makeatletter
\def\@citestyle{\m@th\upshape\mdseries}
\def\citeform#1{{\bfseries#1}}
\def\@cite#1#2{{%
  \@citestyle[\citeform{#1}\if@tempswa, #2\fi]}}
\@ifundefined{cite }{%
  \expandafter\let\csname cite \endcsname\cite
  \edef\cite{\@nx\protect\@xp\@nx\csname cite \endcsname}%
}{}
\makeatother
\begin{document}

\title[MCFs in Warped Products]{Mean Curvature Flows of Closed Hypersurfaces in Warped Product Manifolds}

\author{Zheng Huang, Zhou Zhang, Hengyu Zhou}
\address[Z. ~H.]{Department of Mathematics, The City University of New York, Staten Island, NY 10314, USA}
\address{The Graduate Center, The City University of New York, 365 Fifth Ave., New York, NY 10016, USA}
\email{zheng.huang@csi.cuny.edu}

%\author{Zhou Zhang}
\address[Z. ~Z.]{The School of Mathematics and Statistics, The University of Sydney, NSW 2006, Australia}
\email{zhangou@maths.usyd.edu.au}

%\author{Hengyu Zhou}
\address[H. ~Z.]{College of Mathematics and Statistics, Chongqing University, No. 55, DaxueCheng South Rd., Shapingba, Chongqing, 401331, P. R. China,}
\email{zhouhyu@cqu.edu.cn}
%\author{Zhou Zhang}
%\address[Z. ~Z.]{The School of Mathematics and Statistics, The University of Sydney, NSW 2006, Australia}
%\email{zhangou@maths.usyd.edu.au}


%------------------------------------------------------

\begin{abstract}
We investigate the {\mcf}s in a class of warped products manifolds with closed hypersurfaces fibering over $\R$. In particular, we prove that under some natural conditions on the warping function and Ricci curvature bound for the ambient space, there exists a large class of closed initial hypersurfaces, as geodesic graphs over the {\tg} hypersurface 
$N$, such that the {\mcf} starting from $S_0$ exists for all time and converges to $N$. 
\end{abstract}

\maketitle

%------------------------------------------------------

\section {Introduction}
%%%%%%%%%%%%%%%%%%%
\subsection{Motivation and Main Theorem}
The study of the {\mcf} equation has attracted major attentions in geometric analysis over the past decades. The flow has the following formulation:
\be\label{mcf}	
   \left\{
   \begin{aligned}
      \ppl{}{t}\,F(x,t)&=-H(x,t)\vnu(x,t)\ ,\\
      F(\cdot,0)&=F_{0}\ ,
   \end{aligned}
   \right.
\ene
where $H(x,t)$ and $\vnu(x,t)$ are the {\mc} and unit outward normal vector respectively at $F(x,t)$ of the evolving surface $S(t)$, and our convention of the {\mc} is the sum of the {\pc}s. All other terms will be made transparent later.

A fundamental theorem of Huisken (\cite{Hui84}) in the theory of {\mcf} states that any {\mcf} of a closed and strictly convex initial hypersurface $N^{n-1}\subset \R^n$ (with $n \ge 3$) stays strictly convex, and develops singularity in finite time. This was generalised to several classes of Riemannian manifolds (\cite{Hui86}). Since then there are extensive studies in the field to understand the singularity formulation in various settings and for various {\mcf}s, see for instance \cite{HS99, HS09, LS11, CM12, CIM15} and many others. 

\vspace{0.1in}

We consider a topic similar to that in \cite{Hui84}, namely, under what (natural) conditions, a {\mcf} of closed hypersurfaces in some Riemannian manifolds may converge to some canonical objects. A prototype is a class of {\htm}s (\cite{HLZ16}) where negative curvature and special topology (preventing large balls to appear) help to keep a graphical {\mcf} staying graphical and converging to the {\tg} surface. 

In this paper, we consider a much wider class of {\wppm}s as the ambient space. To describe our setting more precisely, let's first fix notations. Throughout this paper, we use $N$ to denote a closed Riemannian manifold of dimension $n-1$, where $n \ge 3$. We always assume $N$ satisfies the following condition:
\ben
\label{c0}
\item[(C0)] $Ric_N \ge (n-1)\rho g_N$,
\een
for some constant $\rho$. Note that here $\rho$ is different from the one in \cite{Bre13} for simplicity in calculations. The ambient manifold in the current work is $M^n = N^{n-1}\times(-\bar{r},\bar{r})$ for $\bar{r} \in (0,\infty]$, and the {\wpp} metric is 
\be\label{metric}
g = g_M = dr\otimes dr+h^2(r)g_N,
\ene  
with the {\wf} $h(r): (-\bar{r},\bar{r}) \to (0,\infty)$. Geometrically $M$ has the structure of closed hypersurfaces fibering over the real line. The {\wf} $h(r)$ is assumed to satisfy the following conditions:
\ben
\item[(C1)] $h(0) = 1$ and $h'(0) = 0$;
\item[(C2)] $h'(r) > 0$ for all $r\in (0, \bar{r})$ and $h'(r) < 0$ for all $r\in(-\bar{r},0)$;  
\item[(C3)] for $c = \max\{0, \rho\}$ and any $r\in (-\bar{r}, \bar{r})$, 
\be\label{c3}
h(r)h''(r)-h'(r)^2 + \rho \ge c. 
\ene
\een 

The above conditions might seem artificial at the first look, and so let's motivate them below.  

Conditions (C1--2) are natural geometrical conditions. As a consequence, $N$ is {\tg} in $M$, and the unique such hypersurface which is also fixed by the mean curvature flow. The case of $h$ being even provides a natural class of examples in practice. The choice of $h(0)=1$ is just for convenience and certainly not essential.  

Condition (C3) allows negative Ricci curvature on the hypersurface $N$ for suitable {\wf}. Moreover, Condition (C3) ensures the function $|h'(r)/h(r)|$ is non-decreasing in $r$, which is important in the proof of Theorem \ref{graph}. For $r=0$, by (C1), (\ref{c3}) becomes $h''(0)+\rho\ge c$, and so it is easy to generate examples with a proper $\bar r$ value satisfying all the conditions. 

\vspace{0.1in}

Notice that $M$ is not required to be complete as the mean curvature flow in our study stays local in $M$ by the barrier argument as discussed in Lemma \ref{model}. Of course, there are still plenty of examples of complete M, for instance, the one studied in \cite{HLZ16}. Simply speaking, if we choose $N$ to be a closed hyperbolic surface of constant curvature $-1$, then $n = 3$, $h(r) = \cosh(r)$, and $c = \frac12$, then we have $M$ is the {\Fm}, a complete {\htm} as a {\wpp}. All the conditions described above are satisfied. The authors in \cite{HLZ16} used the special structure of the {\Fm} and initiated the study on how graphical {\mcf}s behave in that setting. 

Furthermore, these conditions are comparable with but different from Brendle's in his work on {\cmc} hypersurfaces (\cite{Bre13}), where such interval like $[0, \bar r)$ is considered from general relativity perspective. Indeed, our argument can be applied without any change to the Brendle's setting of $M=N\times [0, \bar r)$ and Conditions (C1--3) will simply be restricted to $[0, \bar r)\subset (-\bar r, \bar r)$. Hence just as in \cite{Bre13}, our main result Theorem \ref{main} can be applied to the de Sitter-Schwarzschild manifold which is of great interest in general relativity. 

\vspace{0.1in}

Let's now fix some notions for the purpose of this paper. 

\vspace{-0.07in}

\begin{Def}\label{angle}
A hypersurface $S$ is called a {\bf graph} over $N$ or {\bf starshaped}, if the {\af} $\Theta = \langle\n,\vnu\rangle > 0$, where $\n = \ppl{}{r}$, and $\vnu$ is the unit normal to $S$. Clearly by definition we have $\Theta \in [0,1]$, and in particular if $\Theta \equiv 1$, we call $S$ is {\bf parallel} (or {\bf equidistant}) to $N$. 
\end{Def}

\vspace{-0.07in}

In this work, we generalize this phenomenon of converging {\mcf} of closed hypersurfaces in \cite{HLZ16} to a much more general class of Riemannian manifolds, namely, those satisfying Conditions (C0--3). More precisely, the following result is proved.

\vspace{-0.05in}

\bt\label{main}
Let $M^n$ be a {\wppm} satisfying Conditions (C0--3) and $S_0$ a smooth closed hypersurface which is a geodesic graph over the unique {\tg} hypersurface $N$ in $M^n$. Then for any $a_0 > 0$, if $S_0$ lies in distance no more than $a_0$ to $N$, and the initial angle satisfies
\be\label{ang-con}
\min_{p\in S_0} \ta(p) \ge  \sqrt{1-\F{1}{h^2(a_0)}},
\ene
 the {\mcf} \eqref{mcf} with the initial hypersurface $S_0$ exists for all time, remains as geodesic graph over $N$ and converges continuously to $N$. Moreover, the convergence is smooth if the above inequality is strict.
\et

\vspace{-0.05in}

Note that the {\wpp} structure in our work is hypersurface fibering over the line. One may also define a type of {\wppm}s as line bundles fibering over hypersurfaces, and some remarkable results of the behaviors of the {\mcf}s in such manifolds were obtained in \cite{BM12}. In that setting, graphs are equidistant graphs, not geodesic graphs. In general we do not expect a geodesic graph stays graphical along the {\mcf}. One can also define the {\mcf} of higher co-dimensions, as well as {\wppm}s where the base manifold is of higher dimensions, and we will not pursue these generalisations here. There are many other interesting papers on various flows in several classes of {\wpp} manifolds (see for example \cite{GLW16, Zho15, Zho16} and others).

%%%%%%%%%%%%%%%%%%%%
%\subsection{Statement of main results}
%------------------------------------------------------

\subsection{Outline of the paper} In Section 2, we discuss important equations and estimates used in the proof of the main result, including the {\ee}s for the height and angle functions along the {\mcf} with a general warped Riemannian manifold as the ambient space. The main result is then proved in Section 3.

%------------------------------------------------------

\subsection{Acknowledgements}  Z. H. thanks a grant from the Simons Foundation (\#359635), and the support from 
U.S. NSF grants DMS 1107452, 1107263, 1107367 ``RNMS: Geometric Structures and Representation varieties" (the GEAR Network), Z. Z. is partially supported by Australian Research Council Future Fellowship FT150100341, and H. Z. is supported by the National Natural Science Foundation of China NSFC \#11801046 and Fundamental Research Funds for the Central Universities, China, Project No. 2019CDXYST0015.

%------------------------------------------------------

\section{Preliminaries}\label{Preliminaries}

In this section, we fix the notations and introduce some preliminary facts that will be used later in this paper.

%------------------------------------------------------

%%%%%%%%%%%%%%%%%%%%%%%%%%%%%%%%%%%%%%%%%%%%

\subsection{Mean curvature flow}\label{mcf-s}

For completeness, we start by collecting and deriving a number of {\ee}s for geometric quantities on $S(t)$, $t\in[0,T)$, which are involved in our calculations. Let $S_0$ be an embedded closed hypersurface in $M^n$, and $F_0$ be the local diffeomorphism representing the embedding: 
\beq
F_0: U \subset \R^{n-1} \to F_0(U) \subset S_0 \subset M.
\eeq

We consider the {\mcf}, i.e. a family of maps $F(\cdot, t)$ satisfying \eqref{mcf}:
\beq
   \left\{
   \begin{aligned}
      \ppl{}{t}\,F(x,t)&=-H(x,t)\vnu(x,t)\ ,\\
      F(\cdot,0)&=F_{0}\, .
   \end{aligned}
   \right.
\eeq
We denote $S(t)$ as the {\ehs} which is the image of the map $F(\cdot, t)$. Here $H(\cdot, t)$ is the {\mc} function of $S(t)$ and $\vnu$ is the correspondingly chosen unit normal of $S(t)$. The short time existence of the {\mcf} was established, and it was further shown that the flow can be extended as long as the norm of the {\sff} is controlled (\cite{Hui84, Hui86}). 

\vspace{0.1in}

As we are primarily working with the graphical case, let's set two functions on $S(t)$: {\it height function} $u(x,t)$, which records the distance, with respect to the metric $g_M$, between any point $x\in S(t)$ and the fixed reference hypersurface $N$,  and the {\it {\af}} (or the gradient function) $\Theta(x,t)$ as defined in Definition {\ref{angle}}. Of course, we always have $\Theta (x,t) \in [0,1]$, and it is clear that the hypersurface $S(t)$ is a geodesic graph over $N$ if $\Theta > 0$ everywhere on $S(t)$ with properly chosen $\vnu$. The general {\ee}s for these two functions are as follows.

\vspace{-0.05in}

\bt\label{old-eq}(\cite {Bar84, EH91, HLZ16}) 
The {\ee}s of $u$ and $\Theta$ have the following form:
\begin{align} \label{eq-u0}
   \ppl{}{t}\,u
        =&\,-H\Theta,\\
               \label{ev-theta}
    \left(\ppl{}{t}-\Delta\right)\Theta=
       &\,(|A|^{2}+ Ric_M(\vnu,\vnu))\Theta+\n(H_{\n})
          -H\inner{\nablabar_{\vnu}\n}{\vnu}\,
\end{align}
where $\n(H_{\n})$ is the variation of {\mc} function of $S(t)$ under the deformation vector field $\n$.
\et

\vspace{-0.05in}

In practice the equation for $\Theta$ is difficult to work with, especially the term $\n(H_{\n})$. In the case of {\Fm}s, a much simplified and explicit equation was derived (\cite{HLZ16}), using special {\hg} of a {\Fm}. In the current situation, we no longer assume the ambient {\wppm} $M$ has constant curvature, and so the {\ee} for $\Theta$ will be more involved. We start with curvature properties in both $M$ and $N$.

%%%%%%%%%%%%%%%%%%%%%
\subsection{{\ricci} in {\wppm}s}

We denote the covariant derivatives of $S(t)$ (induced metric) and $M$ by $\nabla$ and $\nablabar$, respectively. The relationship between the {\ricci}s on $M$ and $N$ is given by the following well-known formula:
\bl\label{ricci}\cite{Bre13}
Let $h=h(r)$ be the {\wf}, then  
\be
Ric_M = Ric_N - (hh''+(n-2)h'^2)g_N - (n-1)\frac{h''}{h}dr\otimes dr.
\ene
\el

\bp This follows from the calculations in \cite{Bes87}. Let $\{e_1, \cdots, e_{n-1}\}$ be a local orthonormal frame on $N$ such that $g_N(e_i,e_j) = \delta_{ij}$. Then one gets
\bear\label{ricci2}
&Ric_M(e_i,e_j) = Ric_N(e_i,e_j) - (hh''+(n-2)h'^2)\delta_{ij}\notag\\
&Ric_M(e_i,\n)=0, \ \ \ Ric_M(\n,\n) = -(n-1)\frac{h''}{h}.
\eear
\ep

Using vectors $\n = \frac{\partial}{\partial r}$ and $\vnu$, it is standard to decompose vector fields into tangential and normal components, {\wrt} either $\n$ or $\vnu$, as follows.
\begin{Def}
For any vector field $X$ in $M$, we define 
\beq
X_\nu = X - \inner X\vnu\vnu,
\eeq
\beq
X_n = X - \inner X\n\n.
\eeq
\end{Def}

With these notations, we have the following decompositions:
\be\label{decom1}
\vnu_n = \vnu - \ta\n,
\ene
\be\label{decom2}
\n_\nu = \n - \ta \vnu = \sum_{k=1}^{n-1}\inner \n{e_k}e_k.
\ene
Clearly, $\vnu_n$ is perpendicular to $\n$ and $\n_\nu$ is perpendicular to $\vnu$. When $\ta = 1$, we have 
$\n = \vnu$ and $\vnu_n = 0$. We further normalize $\vnu_n$ by the metric $g_N$ to set
\be\label{vv}
\vv=\left\{\begin{aligned}
&\F{\vnu_n}{|\vnu_{n}|_{g_N}}\quad &\text{if}\quad \vnu_n\neq 0 \\
&0\quad &\text{if}\quad \vnu_n= 0.
\end{aligned}
\right.
\ene

Now we derive the following technical lemma.

\bl\label{ricci3}
We have 
\be
 Ric_M(\vnu, \n_\nu)=-\F{\ta(1-\ta^{2})}{h^{2}}\{(n-2)(hh''-(h')^{2})+Ric_{N}(\vv, \vv)\}.
\ene
\el
\bp
Let us assume $\vnu_n \neq 0$ or the assertion is trivial. We will omit $M$ in $Ric_M$ for simplicity of notation. By \eqref{decom1} and \eqref{decom2}, we apply \eqref{ricci2} to have
\bear
Ric(\vnu, \n_\nu) &=& Ric(\ta\n+\vnu_n,(1-\ta^2)\n-\ta\vnu_n) \notag\\
&=&\ta(1-\ta^2)Ric(\n,\n) - \ta Ric(\vnu_n,\vnu_n) \notag\\
&=& -(n-1)\ta(1-\ta^2)\frac{h''}{h} - \ta Ric(\vnu_n,\vnu_n).
\eear
Since $\vnu_n$ is tangential to $N$, we find 
\bear
Ric(\vnu_n,\vnu_n) &=& Ric_N(\vnu_n,\vnu_n) - (hh''+(n-2)h'^2)g_N(\vnu_n,\vnu_n)\notag\\
&=& (1-\ta^2)\frac{Ric_{N}(\vv, \vv)}{h^2}-\frac{(1-\ta^2)}{h^2}(hh''+(n-2)h'^2).
\eear
Now the conclusion follows by putting everything above together.
\ep

Another useful fact for our {\wppm} $M$ is that all slices are umbilic. Namely, let $N(a)$ be the equidistant hypersurface  (with signed constant distance $a$) to $N$, and then $N(a)$ is umbilic with constant {\pc} $\frac{h'(a)}{h(a)}$. 
%%%%%%%%%%%%%%%%%%%%%%%%%%%%%%%%%%%%%%%%%%

\subsection{Elliptic equations for height and angle} 

One of the most beautiful geometric properties for {\wppm}s is the existence of the following special vector field which we denote by $V$,  
\be\label{V}
V = h(r)\ppl{}{r} = h(r)\n.
\ene
For any tangential vector field $X$ in $M$, we have (\cite{One83}):
\be\label{XV}
\nablabar_XV = h'(r)X.
\ene
As an application, we immediately have 
\be\label{Xn}
\nablabar_X\n = \frac{h'(r)}{h(r)}(X-\inner X\n\n).
\ene
for any tangential vector field $X$ in $M$. We also calculate the Laplace of the height function $r$ restricted to the {\ehs} $S(t)$ which is denoted by $u$.

\vspace{-0.05in}

 \bpo\label{delta-u}
 Let $\Delta$ be the Laplace operator on the hypersurface $S(t)$. Then we have: 
\be
    \Delta u=\F{h'(u)}{h(u)}(n-2+\ta^{2})-H\ta.
\ene
\epo
\bp
For any point $x \in S=S(t)$, we choose $\{e_1, \cdots,e_{n-1}\}$ (with $e_n= \vnu$) to be a local normal frame of $S$ at $x$. Without loss of generality, we assume that $u(x)\geq 0$.  Then at $x$, we have
\begin{align}
  \Delta u  &=\sum_{i=1}^{n-1} \nb_{e_{i}}\nb_{e_{i}}u \notag\\
                &= \sum_{i=1}^{n-1}  \nb_{e_i} \langle \n, e_i\rangle \notag\\
                &=\sum_{i=1}^{n-1} e_i\inner \n {e_i} \notag\\
                & = \sum_{i=1}^{n-1} \langle \F{h'(u)}{h(u)}(e_i - \langle \n, e_i\rangle \n), e_i \rangle + \sum_{i=1}^{n-1} \langle \n, \bar{\nb}_{e_i} e_i \rangle\notag\\
          & = (n-1)\F{h'(u)}{h(u)} - \F{h'(u)}{h(u)}(1-\Theta^2) - H\Theta  \notag\\
          &= \F{h'(u)}{h(u)}(n-2+\ta^{2})-H\ta,
\end{align}
where we used \eqref{Xn} for the fourth equality.
\ep

An immediate consequence is:

\bcor\label{r}
Using \eqref{eq-u0}, we have the {\ee} for the height function $u$ of $S(t)$ along the {\mcf} :
\be\label{eq-r}
u_t - \Delta u = -\frac{h'}{h}(n-2+\ta^2).
\ene
\ecor

We now derive the main technical tool in this work, the calculation for the Laplacian of the {\af} $\ta$.

\bt\label{Dta}
We have
\begin{align*}
   \Delta \ta&= \inner{\nabla H}\n  -|A|^{2}\ta +\F{h'(r)}{h(r)}\{H(\ta^{2}+1)-2\inner{\n}{\nabla\ta}\} -(n-1)\F{h'(r)^2}{h^2(r)}\ta\notag\\
   &\ \  - \F{\ta(1-\ta^{2})}{h^{2}(r)}[(n-1)(h(r)h''(r)-h'(r)^2)+Ric_{N}(\vv, \vv)],
\end{align*}
where $\vv$ is defined in \eqref{vv}.
\et

\bp
We first work with the auxilliary function:
\be\label{eta}
\eta = \inner{V}\vnu = h(r)\ta.
\ene
We still use the local normal frame $\{e_1, \cdots,e_{n-1}\}$ (with $e_n= \vnu$) at any point $x \in S=S(t)$ such that at this point $x$, we have 
\beq
\nabla_{e_i}e_k(x) = 0 (i\neq k), \ \ \ \bnb_{e_i}e_i(x) = -a_{ii}\vnu.
\eeq
where $A = (a_{ij})$ is the {\sff} of the hypersurface $S$. 

We then have the following computation at $x$:
\begin{align*}
     \Delta \eta &= \sum_{i=1}^{n-1} \nb_{e_{i}}\nb_{e_{i}}\inner V\vnu \notag\\
&= \sum_{i=1}^{n-1} \inner{\bnb_{e_{i}}\bnb_{e_{i}}V}\vnu+ 2\inner{\bnb_{e_{i}}V}{\bnb_{e_{i}}\vnu}+ \inner V{\bnb_{e_{i}}\bnb_{e_{i}}\vnu}\notag\\
            & = \sum_{i=1}^{n-1} \inner{\bnb_{e_{i}}(h'(u)e_{i})}\vnu + 2 \sum_{i,k=1}^{n-1}h'(u)\inner{e_{i}}{a_{ik}e_k}
            +\sum_{i,k=1}^{n-1} \inner V{\bnb_{e_i}(a_{ik}e_k)}  \notag\\
            & = -h'(u)H + 2h'(u)H + \sum_{i,k=1}^{n-1} a_{ik}\inner V{\bnb_{e_i}e_k}+a_{ik, i}\inner V{e_k} \notag\\
            & =h'(u)H -|A|^2\eta+ \sum_{i,k=1}^{n-1}a_{ik, i}\inner V{e_k}\,
\end{align*}
where we have made use of the properties of the normal frame at the point $x$ under investigation. Let's examine the last summation more closely. First we use \eqref{decom2} to get
\beq
     \sum_{k=1}^{n-1}\inner V{e_k}e_k  =h(r) \sum_{k=1}^{n-1}\inner \n{e_k}e_k  = h(r)\n_\nu.
\eeq
Then recall the Codazzi equation for $S \subset M$:
\be\label{codazzi}
e_i(a_{ik}) - e_k(a_{ii})=\bar{R}(\vnu,e_i,e_k,e_i),
\ene
where $\bar{R}$ is the curvature tensor on $M$. Thus, we have
     \begin{align*}
     \sum_{i,k=1}^{n-1}a_{ik, i}\inner V{e_k} & =\sum_{i,k=1}^{n-1}e_k(a_{ii})\inner V{e_k} + \bar{R}(\vnu,e_i,e_k,e_i)\inner V{e_k} \notag\\
&= \inner V{\nabla H}+ \sum_{i,k=1}^{n-1}\bar{R}(\vnu,e_i, \inner V{e_k}e_k,e_i) \notag\\
            & = \inner V{\nabla H} + h(r)\sum_{i=1}^{n-1} \bar{R}(\vnu,e_i, \n_\nu,e_i)\notag\\
            & =\inner V{\nabla H} + h(r)Ric_M(\vnu, \n_\nu)\,.
     \end{align*}
Now we are in position to apply Lemma ~\ref{ricci3} and arrive at:
\begin{align} \label{etau}
   \Delta \eta &= h'(r)H -|A|^{2}\eta + \inner V{\nb H}\notag\\
   &\ \ \ \ \ -\F{\ta(1-\ta^{2})}{h(r)}[(n-2)(h(r)h''(r)-(h'(r))^{2})+Ric_{N}(\vv, \vv)].
\end{align}
Since $\eta = h(r)\ta$, and applying Lemma ~\ref{delta-u}, we have 
\begin{align*}
\Delta \eta & =h\Delta\ta+ 2h'\inner\n{\nabla\ta}+\ta(h'\Delta r+h''|\nabla r|^2)  \notag\\
&=h\Delta\ta+ 2h'\inner\n{\nabla\ta}+\ta(h'\Delta r+h''(1-\ta^2))  \notag\\
&=h\Delta\ta+ 2h'\inner\n{\nabla\ta}+\ta[\F{h'^2}{h}(n-2+\ta^2)-h'H\ta+h''(1-\ta^2)]  \notag\,.
\end{align*}
Now the assertion follows by isolating $\Delta\ta$ on one side of the equation.
\ep
%%%%%%%%%%%%%%%%%%%%%%%%%%%%%%
\section{Proof of Main Theorem}

The key is to obtain a positive lower bound for the angle function $\Theta$. 

%%%%%%%%%%%%%%%%%%%%%%%%%%%%%%%
\subsection{Evolution equation for the angle function squared} 

In Theorem ~\ref{Dta}, we have derived the Laplacian of $\ta$. Now let us derive the {\ee}.

\vspace{-0.05in}

\bt\label{ee-ta}
The {\af} $\ta(\cdot,t)$ satisfies the following {\ee} along the {\mcf} \eqref{mcf}:
\begin{align}
   \ta_t- \Delta \ta&= |A|^{2}\ta +\F{2h'(r)}{h(r)}\{\inner{\n}{\nabla\ta} - H\} +(n-1)\F{h'(r)^2}{h^2(r)}\ta\notag\\
   &\ \ \ +\F{\ta(1-\ta^{2})}{h^{2}(r)}\{(n-1)(h(r)h''(r)-h'(r)^2)+Ric_{N}(\vv, \vv)\},
\end{align}
where $\vv$ is defined in \eqref{vv}.
\et

\bp
For the {\mcf} we have (\cite{Hui86}) $\ppl{}{t}\,\vnu=\nabla H$. Then we have 
\bear
   \ppl{\Theta}{t} &=& \ppl{}{t}\inner\n\vnu =\inner{\ppl{}{t}\vnu}\n + \inner\vnu{\bnb_{-H\vnu}\n}  \notag\\
   &=& \inner{\nb H}\n -H\inner\vnu{\bnb_{\vnu}\n}.\notag
\eear
Using \eqref{Xn}, we have
\beq
\bnb_{\vnu}\n = \F{h'}{h}(\vnu - \Theta\n) =  \F{h'}{h}\vnu_n,
\eeq
and so that
\be\label{theta-t}
\ppl{\Theta}{t} = \langle\nb H, \n\rangle -\F{Hh'}{h}(1-\Theta^2).
\ene
Now the assertion follows from combining this with Theorem ~\ref{Dta}.
\ep

It's actually easier to work with $\Theta^2$. Let's denote $f(\cdot,t) = \Theta^2$ which satisfies the following differential inequality.
\bcor\label{f-ineq}
As long as $f(\cdot,t) = \Theta^2 > 0$ and the inequality \eqref{c3} holds, we have 
\be\label{ineq-f}
   f_t- \Delta f \ge \inner{\F{2h'(r)}{h(r)}\n-\F{\nabla f}{2f}}{\nabla f} +G(f,r),
\ene
where
\be\label{G}
G(f,r) = \F{2(n-1)(1-f)}{h^2(r)}\{cf-h'(r)^2\},
\ene
and $c$ is the constant from condition (C3).
\ecor
\bp
A direct calculation from Theorem ~\ref{ee-ta} shows that $f$ satisfies the following {\ee}:
\begin{align}\label{ta-square}
   f_t- \Delta f&= 2|A|^{2}f +\inner{\F{2h'(r)}{h(r)}\n-\F{\nabla f}{2f}}{\nabla f} - \F{4h'(r)\sqrt{f}}{h(r)}H+2(n-1)\F{h'(r)^2}{h^2(r)}f\notag\\
   &\ \ \ +\F{2f(1-f)}{h^{2}(r)}\{(n-1)(h(r)h''(r)-h'(r)^2)+Ric_{N}(\vv, \vv)\}.
\end{align}
By \eqref{c3}, we have 
\begin{align}
   f_t- \Delta f&\ge 2|A|^{2}f +\inner{\F{2h'(r)}{h(r)}\n-\F{\nabla f}{2f}}{\nabla f} - \F{4h'(r)\sqrt{f}}{h(r)}H+2(n-1)\F{h'(r)^2}{h^2(r)}f\notag\\
&\ \ \ +\F{2f(1-f)}{h^{2}(r)}\{(n-1)c\},
\end{align}
where the special case of $\vv=0$ is taken care of by the factor $1-f$ in the last term on the right hand side of (\ref{ta-square}) because $f=\Theta^2=1$ in this case.  

Now the corollary follows from using the fact that $|H|\le \sqrt{n-1}|A|$ and completing the square. 
\ep

%%%%%%%%%%%%%%%%%%%%%%%%%%
\subsection{Barriers}

In this subsection, we treat the model case, namely, the {\mcf} of {\ihs} parallel (i.e equidistant) to $N$. By the well-known {\ap} in {\mcf}s of closed hypersurfaces, this special {\mcf} will serve as barriers to control the behavior of our {\mcf} for a more general class of initial hypersurfaces.

\vspace{-0.05in}

\bl\label{model}
Let $N$ be a closed $(n-1)$-dimensional Riemannian manifold and $M = N\times \R~\text{or}~N\times [0, \infty)$ the {\wppm} with the metric given by \eqref{metric} satisfying Conditions (C1--2). Then any {\mcf} \eqref{mcf} in $M$ of {\ihs} $N(a)$, where $N(a)$ is a hypersurface of constant (signed) distance $a \in \R$ to $N$, exists for all time, stays umbilic and converges smoothly to $N$ as $t \to \infty$.
\el

\bp
Since $N(a)$ is parallel to $N$, the initial {\af} $\ta(0)\equiv 1$, moreover, $N(a)$ is umbilic with {\pc} $\F{h'(a)}{h(a)}$. By uniqueness of the {\mcf}, we have $\ta \equiv 1$ and the evolving surface stays parallel to $N$.

It's trivial when $a=0$. Let's assume $a >0$. Let $R(t)$ be the height function of the {\ehs} at time $t$. Note that it is a function of $t$ only, since $N(R)$ is parallel to $N$. Then by either \eqref{eq-u0} or \eqref{eq-r}, we have 
\be\label{ode4u}
  \left\{
   \begin{aligned}
      \ddl{R(t)}{t} = -(n-1)\F{h'(R(t))}{h(R(t))}\ \\
      R(0) =a >0\ .
   \end{aligned}
   \right.
\ene
Since $h$ is a positive function over $\mathbb{R}$, the solution $R(t)$ exists forever. In light of the direction field of this ODE, mostly just the sign of derivative, we know that $R(t)$ decreases and stays positive. Assuming $R(t)\to A\ge 0$ as $t\to \infty$, one can easily rule out the case of $A>0$ since $\frac{h'(A)}{h(A)}>0$, so $A=0$ and $N(R)$ converges smoothly to $N$ at time infinity. The case of $a<0$ can be treated in the same way. In the case of $M=N\times [0, \infty)$, $N=N\times\{0\}$ itself serves as the other barrier. This completes the proof. 
\ep

As a corollary, we obtain the convergence part of the main theorem. 
\bcor\label{conv}
Let \eqref{mcf} be any {\mcf} in the above {\wppm} $M$ with a closed {\ihs} $S(0)$. If it exists for all time, then it converges continuously to $N$.
\ecor
\bp
Since the {\ihs} is closed, there is a constant $a > 0$ such that $S(0)$ is enclosed in the region between parallel hypersurfaces $N(-a)$ and $N(a)$. The corollary follows from the fact that {\mcf}s of {\ihs}s $N(-a)$ and $N(a)$ both converge to $N$ and the {\ap}. 
\ep

We conclude this subsection by the following lemma which will be used later to obtain the key estimate.

\bl\label{con}
Let $R(t)> 0$ be the solution for the {\ivp} \eqref{ode4u}, and $\bar{f}(t)$ be the solution to the {\ivp}:
\be\label{ode4f}
  \left\{
   \begin{aligned}
      \ddl{\bar{f}(t)}{t} = -2(n-1)(1-\bar{f})\F{h'(R(t))^2}{h^2(R(t))}\ \\
      \bar{f}(0) = {\bar{f}}_0  \in [0,1]\ ,
   \end{aligned}
   \right.
\ene
Then we have the following identity:
\be\label{iden}
(1-\bar{f}(t))h^2(R(t)) = (1-\bar{f}_0)h^2(a),
\ene
for all $t \ge 0$.
\el

\bp
Let's assume $\bar{f}_0 < 1$, otherwise the equation \eqref{ode4f} forces $\bar{f}(t) = 1$ and we are done. Since we have both $h(r) > 0$ and $h'(r) > 0$ for all $r > 0$, then 
\beq 
\ddl{\bar{f}(t)}{t} = -2(n-1)(1-\bar{f})\F{h'^2(R(t))}{h^2(R(t))} < 0.
\eeq
So the solution to \eqref{ode4f} exists for all time. To prove the lemma, we set the function $\Lambda(t) = (1-\bar{f}(t))h^2(R(t))$ and prove that it is actually independent of $t$. We justify this by a direct calculation using both equations \eqref{ode4u} and \eqref{ode4f}:
\bear
\ddl{\Lambda(t)}{t} &=& -h^2(R(t))\ddl{\bar{f}(t)}{t} + (1-\bar{f}(t))(2h(R(t))h'(R(t))\ddl{R(t)}{t}) \notag\\
&=& 2(n-1)(1-\bar{f})h'^2(R(t)) +2(1-\bar{f})hh'(-(n-1)\F{h'(R(t))}{h(R(t))})  \notag\\
&=& 0.
\eear
\ep

\br \label{lim} As an immediate consequence, we have the limit 
\be
\lim_{t\to\infty}\bar{f}(t) = 1 - (1-\bar{f}_0)h^2(a).
\ene
This is the only place that $h(0)=1$ is used which can of course be easily adjusted using any positive constant instead.  
\er

%%%%%%%%%%%%%%%%%%%%%%%%%%%%%%%%%%%%
\subsection{Gradient estimate} 

Now we apply the barriers established in the previous subsection and comparison equations to control the lower bound for $\ta$ and establish the gradient estimate. In particular, we prove: 

\vspace{-0.05in}

\bt\label{graph}
Let $M^n$ be a {\wppm} satisfying Conditions (C0--3) and $S_0$ a smooth closed hypersurface which is a geodesic graph over the unique {\tg} hypersurface $\Sigma$ in $M^n$, and suppose there is a constant $a_0 > 0$ such that $S_0$ lies between $\Sigma(\pm a_0)$. Then if $S_0$ satisfies 
$$\min_{p\in S_0} \ta(p) \ge \sqrt{1-\F{1}{h^2(a_0)}},$$ 
the {\mcf} \eqref{mcf} with initial hypersurface $S_0$ remains graphical over $\Sigma$, namely, $\ta(\cdot, t) > 0$ as long as the flow exists.
\et

\bp
Let us first recall the evolution inequality satisfied by $f = \ta^2$ as in Corollary ~\ref{f-ineq}:
\beq
   f_t- \Delta f \ge \inner{\F{2h'(r)}{h(r)}\n-\F{\nabla f}{2f}}{\nabla f} +G(f,r),
\eeq
where
\beq
G(f,r) = \F{2(n-1)(1-f)}{h^2(r)}\{cf-h'(r)^2\},
\eeq
and $c$ is the constant from condition (C3). Let $\phi(t)$ be the spatial minimum of $f(\cdot, t)$ on the {\ehs} $S(t)$, namely, 
$$\phi(t)=\min_{S(t)} f. $$ We only have to establish $\phi \in (0,1]$ for a priori estimate. At the spatial minimum of $f$, we have 
$\nabla f=0$ and $\Delta f\ge 0$, and so for $t>0$ (using Hamilton's trick), we find:
\begin{align}
 \F{d\phi}{dt} & \ge \ppl{f}{t}  - \Delta f \notag\\
     &\ge \F{2(n-1)(1-f)}{h^2(r)}\{cf-h'(r)^2\} \notag\\
     &\ge -\F{2(n-1)(1-f)}{h^2(r)}\{h'(r)^2\} \notag \\
     &= -2(n-1)(1-\phi)\big\{\F{h'(r)}{h(r)}\big\}^2. \notag
\end{align}
By our conditions on the warping function $h(r)$, in particular \eqref{c3}, we have $\vline\F{h'(r)}{h(r)}\vline$ is nondecreasing in $|r|$. Let $R(t) > 0$ be the solution to the {\ivp} \eqref{ode4u}, namely, the evolving distance of the {\mcf} with {\ihs} $N(a_0)$ at any time $t \ge 0$. Since our {\ehs}s $S(t)$ in the {\mcf} \eqref{mcf} of {\ihs} $S_0$ is squeezed by the barriers $N(R(t))$, by Lemma ~\ref{model}, we have $|r(t)| \le R(t)$, and therefore 
\be\label{phi-ineq}
 \left\{
   \begin{aligned}
     \F{1}{1-\phi} \F{d\phi}{dt} \ge -2(n-1)\F{h'^2(R(t))}{h^2(R(t))}\ \\
      \phi(0) = \phi_0  \in (0,1)\ .
   \end{aligned}
   \right.
\ene
Recall that our comparison {\ivp} \eqref{ode4f} is equivalent to:  
\beq
  \left\{
   \begin{aligned}
      \F{1}{1-\bar{f}} \ddl{\bar{f}(t)}{t} = -2(n-1)\F{h'^2(R(t))}{h^2(R(t))}\ \\
      \bar{f}(0) = {\bar{f}}_0  \in [0,1]\ ,
   \end{aligned}
   \right.
\eeq
Choosing $\phi_0=\bar f_0$, we find, for $t > 0$, 
$$\frac{d}{dt}\left(\log\F{1-\bar{f}(t)}{1-\phi(t)}\right) \ge 0,$$ 
which implies $\phi(t) \ge \bar{f}(t)$. Finally, by Lemma ~\ref{con} and Remark ~\ref{lim}, setting $a = a_0$, we have 
$$\phi(t) \ge \bar{f}(t) \ge\lim_{t\to\infty}\bar{f}(t) = 1 - (1-\bar{f}_0)h^2(a_0).$$
where the second $\ge$ is strict unless the initial surface is $N$, for which there is nothing to proof. 

Now we apply the initial angle condition $\min_{p\in S_0} \ta(p) \ge \sqrt{1-\F{1}{h^2(a_0)}}$, and choose $\bar{f}_0$ to be $(\min_{p\in S_0} \ta(p))^2$, we arrive at $\phi(t) > 0$ as long as the flows exists. This completes the proof.
\ep

%%%%%%%%%%%%%%%%%%%%%%
\subsection{Completing the proof of the main result}

Now we can assemble the ingredients and complete the proof of Theorem \ref{main}.

\bp (of Theorem ~\ref{main})
We have shown in Theorem ~\ref{graph} that the {\mcf} \eqref{mcf} stays graphical as long as it exists. This provides the gradient estimate for the {\mcf} for any finite time interval. By the classical theory of parabolic equations in divergent form (for instance \cite{LSU68}), the higher regularity and a priori estimates of the solution follow in the standard way. This yields the long time existence of the flow by Huisken (\cite{Hui86}). Thus, by Lemma \ref{model} and the avoidance principle, the continuous convergence of the flow also follows. 

When the inequality in the assumption of the theorem is strict, the proof of Theorem ~\ref{graph} gives a uniform positive lower bound of the angle for all time, and so the higher order estimates are uniform for all time, providing the smooth convergence. Hence, the proof of Theorem ~\ref{main} is completed.
\ep

\subsection{Applications and remarks}

We begin by showing that the main result can be applied for the de Sitter-Schwarzschild manifold from general relativity. The discussion is adjusted from the same consideration in \cite{Bre13}. Usually in literature, such manifold is described in the following setting: $M=S^{n-1}\times (\underline s, \overline s)$ with the warped metric 
$$g=\frac{1}{\omega(s)}ds\otimes ds+s^2 g_{_{S^{n-1}}}$$
where $S^{n-1}$ is the unit $(n-1)$-sphere, $\omega(s)>0$ in $(\underline s, \overline s)\subset (0, \infty)$ and smooth up to $\underline s$ with $\omega(\underline s)=0$. $Ric_{S^{n-1}}= (n-2)g_{_{S^{n-1}}}$ gives Condition (C0) with $\rho=\frac{n-2}{n-1}$. For the de Sitter-Schwarzschild manifold, 
$$\omega(s)=1-ms^{2-n}-\kappa s^2$$  
where the mass constant $m>0$, the cosmological constant $\kappa$ can be either non-positive in which case $\overline s=\infty$ or satisfy
$$n^nm^2\kappa^{n-2}<4(n-2)^{n-2}$$ 
so that $\omega(\overline s)=0$ for $0<\underline s<\overline s<\infty$. We then use the change of variable in \cite{Bre13}, $r=F(s)$ with 
$$\frac{dr}{ds}=F'(s)=\frac{1}{\sqrt{\omega(s)}}, ~~F(\underline s)=0.$$
The assumption on $\omega(s)$ makes sure that this is a legitimate change of variable from $s\in [\underline s, \overline s)$ to $r\in [0, F(\overline s))$. This provides a natural way to extend the manifold to $s=\underline s$ or $r=0$. 

Using the new warping variable $r$, we have  
$$g=dr\otimes dr+h^2(r)g_{_{S^{n-1}}},$$
where $h(r)=s=F^{-1}(r)$. So we have 
$$h'(r)=\frac{ds}{dr}=\sqrt{\omega(s)}$$
which clearly satisfies Conditions (C1--2) except that $h(0)=\underline s$ which is not a problem as described before. Moreover, we have
$$h''(r)=\frac{1}{2\sqrt{\omega(s)}}\omega'(s)\frac{ds}{dr}=\frac{\omega'(s)}{2}.$$
Regarding Condition (C3), $c=\rho>0$, we need
$$h(r)h''(r)-h'^2(r)=\frac{1}{2}s~\omega'(s)-\omega(s)\geq 0$$
and so at $r=0$, i.e. $s=\underline s$, 
$$h(0)h''(0)-h'^2(0)=\frac{1}{2}\underline s~\omega'(\underline s).$$
For the de Sitter-Schwarzschild manifold, 
$$\omega'(s)=-m(2-n)s^{1-n}-2\kappa s,$$
and so we have
$$h(r)h''(r)-h'^2(r)=\frac{1}{2}s~\omega'(s)-\omega(s)=\frac{1}{2}mns^{2-n}-1.$$
Hence we need $\frac{1}{2}mns^{2-n}-1\geq 0$, and so 
$$s\leq \left(\frac{mn}{2}\right)^{\frac{1}{n-2}}.$$

\noindent {\it Claim: $\left(\frac{mn}{2}\right)^{\frac{1}{n-2}}\in (\underline s, \overline s)$.} 

\begin{proof} There are two cases to deal with.  

If $\kappa\leq 0$, since $N\geq 3$, it's obvious that 
$$h(0)h''(0)-h'^2(0)=\frac{1}{2}\underline s~\omega'(\underline s)> 0.$$
So $\left(\frac{mn}{2}\right)^{\frac{1}{n-2}}\in (\underline s, \infty)=(\underline s, \overline s)$. 

\vspace{0.1in}

If $\kappa>0$ and $n^nm^2\kappa^{n-2}<4(n-2)^{n-2}$, as $s$ increases from $0$ to $\infty$, it's clear that $\omega'(s)=-m(2-n)s^{1-n}-2\kappa s$ strictly decreases from $\infty$ to $-\infty$, and so we know that $\omega(s)$ strictly increases (passing $\omega(\underline s)=0$) from $-\infty$ to its positive maximum which is achieves at some $S\in (\underline s, \overline s)$, and then strictly decreases (passing $\omega(\overline s)=0$) to $-\infty$. So we clearly have $\omega'(\underline s)>0$, and  
$$h(0)h''(0)-h'^2(0)=\frac{1}{2}\underline s~\omega'(\underline s)> 0.$$
So $\left(\frac{mn}{2}\right)^{\frac{1}{n-2}}>\underline s$. Meanwhile, at $s=S$ where the maximum of $\omega(s)$ is achieved, 
$$\frac{1}{2}S~\omega'(S)-\omega(S)=-\omega(S)<0,$$
and so $\overline s>S>\left(\frac{mn}{2}\right)^{\frac{1}{n-2}}$.

So the claim is justified. 

\end{proof}

Now we can conclude that Condition (C3) is satisfied in $S^{n-1}\times [\underline s, \left(\frac{mn}{2}\right)^{\frac{1}{n-2}}]$ as part of the de Sitter-Schwarzschild manifold, where we can apply the main result of Theorem \ref{main}.  

\vspace{0.1in}

Finally, we would like to point out that there are examples where graphical complete hypersurfaces of {\wppm}s fail to stay graphical along the {\mcf} after some finite time, for example, as discussed in Appendix A of \cite{Zho17}. A priori it's not clear whether this will indicate the development of geometric singularities along the flow. We hope to address this problem in future works.   

%%%%%%%%%%%%%%%%%%%%

\bibliographystyle{amsalpha}
\bibliography{ref-warp}

\end{document}